\documentclass[preprint,12pt]{elsarticle}

\newtheorem{theorem}{Theorem}[section]

\newtheorem{corollary}{Corollary}[section]
\newtheorem{proposition}{Proposition}[section]

\newtheorem{remark}{Remark}[section]
\newcommand{\ignore}[1]{}{}

\usepackage{amsmath}
\usepackage{amssymb}

\usepackage{color}
\usepackage{soul}
\usepackage[dvipsnames]{xcolor}

\usepackage[colorlinks=true, linkcolor=blue, citecolor=blue]{hyperref}

\def\1{{{\mbox{${\rm{1\negthinspace\negthinspace I}}$}}}}

\newcommand\beq{\begin{equation}}
\newcommand\eeq{\end{equation}}

\usepackage{ifthen}
\usepackage{xkeyval}
\usepackage{todonotes}
\begin{document}

\begin{frontmatter}

\title{Gaussian martingale inequality applies to random functions and  maxima of empirical processes}
\author{Xiequan Fan}
 \cortext[cor1]{\noindent Corresponding author. \\
\mbox{\ \ \ \ }\textit{E-mail}: fanxiequan@hotmail.com (X. Fan). }
\address[cor1]{Center for Applied Mathematics, Tianjin University,  300072 Tianjin, China }
\address[cor2]{Regularity Team, Inria, France}

\begin{abstract}
We obtain  a Bernstein type Gaussian concentration inequality for martingales.
 Our inequality improves  Azuma-Hoeffding's inequality for moderate deviations $x$.
Following the work of   \mbox{McDiarmid} \cite{M89}, Talagrand \cite{T96} and Boucheron,   \mbox{Lugosi} and Massart \cite{BLM00,BLM03}, we show that our result can be applied to the concentration of random functions,  Erd\"{o}s-R\'{e}nyi  random graph, and maxima of empirical processes. Several interesting Gaussian concentration inequalities    have been obtained.
\end{abstract}

\begin{keyword} martingales; concentration inequalities;  Azuma-Hoeffding's inequality; McDiarmid's inequality;   Freedman's inequality; Self-bounding functions; Rademacher averages; maxima of empirical processes
\vspace{0.3cm}
\MSC Primary 60G42; 60E15; 60F10;  Secondary 62G30
\end{keyword}

\end{frontmatter}




\section{Introduction}
Concentration inequalities for tail probabilities of  random variables play an important role both in theoretical study
and applications. Several methods have been introduced
for obtaining such inequalities, including Bernstein's method \cite{B46} (see also   Hoeffding \cite{Ho63}),  Talagrand's induction method \cite{T96,T97},  and the powerful so-called ``entropy method", based on logarithmic Sobolev inequalities, developed by  Ledoux \cite{L96}.
The  martingale method is also a useful way to derive  concentration inequalities; see McDiarmid \cite{M89,M98} for an excellent survey
and also Bercu, Delyon and Rio \cite{BDR08} for a recent monograph. For instance, when the martingale differences are bounded,   Azuma-Hoeffding's inequality is one of the classical results.

Assume that we are given a finite sequence of centered real-valued random variables $(X _i)_{i=1,...,n}.$ Let $S_k=\sum_{i=1}^kX_i$ be the partial sums of $(X _i)_{i=1,...,n}.$  If $(X _i)_{i=1,...,n}$ are independent and
satisfy $|X _i|\leq a_i$ for some constants  $a_i,$ $i=1,...,n,$   Azuma-Hoeffding's inequality \cite{Ho63}  implies that,
for any positive $x,$
\begin{eqnarray}\label{ineq3}
 \mathbf{P}\Big(S_n  \geq x  \Big)
 \leq  \exp \left\{-\frac{x^2}{ 2 \sum_{i=1}^n a_i^2} \right\} .
\end{eqnarray}
Denote by $\textrm{Var}(X)$ the variance of  $X.$ Notice that the factor $\sum_{i=1}^na_i^2$  is an upper bound of
the variance of  $S_n$, that is $\textrm{Var}(S_n) \leq \sum_{i=1}^n a_i^2$. Since $\textrm{Var}(S_n)=\sum_{i=1}^n \textrm{Var}(X_i)$,
it is easy to see that  Azuma-Hoeffding's inequality is Gaussian provided that $\textrm{Var}(X_i)=a_i^2$ for all $i=1,...,n$.
Otherwise, $\textrm{Var}(S_n) <\sum_{i=1}^n a_i^2$ which means   Azuma-Hoeffding's inequality is sub-Gaussian and not tight enough.
 The following Bernstein inequality  is significantly better than  Azuma-Hoeffding's inequality.  Bernstein \cite{B46}  proved that, for any positive $x,$
\begin{eqnarray}\label{ineq4}
 \mathbf{P}\Big(S_n  \geq x  \Big)
 \leq  \exp \left\{-\frac{x^2}{ 2 (\textrm{Var}(S_n) + \frac13 a \, x )} \right\} ,
\end{eqnarray}
where $a= \max\{a_i, i=1,...,n\}.$
It is obvious that Bernstein's bound is smaller than Azuma-Hoeffding's bound for all $0\leq x \leq \frac3a\big(\sum_{i=1}^n a_i^2-\textrm{Var}(S_n)\big).$
Moreover, Bernstein's bound (\ref{ineq4}) is a Gaussian  bound for moderate deviations $x$, in contrast to  (\ref{ineq3}) which is sub-Gaussian.   That is why   Bernstein's inequality is significantly better than  Azuma-Hoeffding's inequality for moderate deviations $x$.

In the sequel, consider the martingale case.  Assume that $(X _i,\mathcal{F}_i)_{i=1,...,n}$ is a sequence of real-valued martingale differences,
which means $(S_k,\mathcal{F}_k)_{k=1,...,n}$ is a  martingale. Assume  $|X_i| \leq a_i$ for all $i=1,...,n.$   Inequality (\ref{ineq3}) also holds true; see \cite{A67,Ho63}. This martingale version of   Azuma-Hoeffding's inequality   has  many interesting applications;  see \mbox{McDiarmid} \cite{M89,M98} for concentration of functions of independent random variables   and see Devroye \cite{D91} for   nonparametric estimations.
Notice that for bounded martingale differences $|X_i| \leq a_i$ for all $i=1,...,n$, the term $\sum_{i=1}^na_i^2$  is also the upper bound of
the variance of $S_n$.   Thus  we also would like to replace it by a smaller one. Freedman  \cite{FR75} (see also \mbox{van de Geer} \cite{V95}, de la Pe\~{n}a \cite{D99} and   \cite{F12,FGL15} for closely related results) has obtained
the following inequality for martingales: for any positive $x,$
\begin{eqnarray}
 \mathbf{P}\Big(  S_n  \geq x\ \Big)
  \leq   \exp \left\{-\frac{x^2}{ 2(||\!\left\langle S\right\rangle_n\!||_{\infty}+ \frac13 a\, x  )} \right\} , \label{freedmanxf2}
\end{eqnarray}
where $a= \max\{a_i, i=1,...,n\} $ and $\left\langle S\right\rangle_n=\sum_{i=1}^n\mathbf{E}[X_{i}^2|\mathcal{F}_{i-1}]$ is the predictable quadratic variance of $S_n.$  Notice that $\left\langle S\right\rangle_n$ can be substantially
smaller that $\sum_{i=1}^n a_i^2$. Thus   Freedman's  inequality provides a sharper bound than   Azuma-Hoeffding's inequality for moderate deviations $x$.  However,    Freedman's  inequality does not generality provide a  Gaussian   bound.
The reason is that  $||\left\langle S\right\rangle_n||_{\infty}$ is usually larger than the variance of $S_n$ and   Freedman's  inequality can be much worse than its ``Gaussian'' version which should make $\textrm{Var}(S_n)$ appear instead of
$||\left\langle S\right\rangle_n||_{\infty}$.

In this paper  we would like to establish 
 a  ``Gaussian''  version of   Bernstein's inequality for martingales.   We prove that  if the martingale differences satisfy the following Bernstein type condition
 \begin{eqnarray} \label{c0}
 \mathbf{E}[|X_i|^k ] \leq \frac{1}{2} \frac{  k! \, a ^{k-2}  }{(k-1)^{k/2}}  \mathbf{E}[X_i^2],\ \ \ \  \ k\geq  2,
\end{eqnarray}
for some positive constant $a,$ then, by  inequality (\ref{jnsk}), for any positive $x,$
\begin{eqnarray}\label{fsf}
 \mathbf{P}\Big(  S_n  \geq x  \Big)
  \leq  \exp \left\{-\frac{x^2}{ 2(  \textrm{Var}(S_n) + a x \sqrt{n} \, )} \right\} .
\end{eqnarray}
We will show that both the bounded difference $|X_i|\leq a$ and the normal random random variables  satisfy condition (\ref{c0}). It is interesting to see that our inequality (\ref{fsf}) is a  Gaussian   bound. Thus (\ref{fsf}) has certain advantages over  Azuma-Hoeffding's inequality.
Moreover, we also prove that the constant $a$ in (\ref{fsf}) cannot be replaced by  positive constants $b_n,$ with $b_n\rightarrow 0$ as $n\rightarrow \infty$, under the stated condition. In this sense,  inequality (\ref{c0}) can be regarded as  a true martingale  version of  Bernstein's  inequality.

Assume that $\textrm{Var}(S_n)/n$ is bounded. Inequality (\ref{fsf}) implies that
\begin{eqnarray}
\ \ \mathbf{P}\Big(    S_n  \geq  x\sqrt{n} \Big)  =  O\Big( \exp\{ - C_1 \, x\} \Big),\ \ \ \ \ \ x\rightarrow \infty,
\end{eqnarray}
for some positive constant $C_1$ which does not depend on $x$. The last equality is also the best possible in the sense that there exists a class of
stationary  martingale differences such that (\ref{c0}) holds and
\begin{eqnarray}
\mathbf{P}\Big(    S_n  \geq  x\sqrt{n} \Big)   \geq \exp\Big\{ - C_2\, x \Big\},\ \ \ \ \ \ x\rightarrow \infty,
\end{eqnarray}
for some positive constant $C_2$ which does not depend on $x$.

Inequality (\ref{fsf}) is user-friendly.
Following the work of McDiarmid  \cite{M89}, we apply inequality (\ref{fsf}) to Lipshitz functions of independent random variables.
Following the work of Boucheron, Lugosi and Massart \cite{BLM00,BLM03,B09}, and \mbox{McDiarmid} and Reed  \cite{MR0}, we apply our
inequalities to self-bounding functions of independent random variables, Vapnik-Chervonenkis entropies, Rademacher averages,  and counting small subgraphs in random graphs.
 At last but not least, applying our inequalities to the maxima of empirical processes,
we provide a partial positive  answer to Massart's question \cite{M00} about the best constants in   Talagrand's  inequality \cite{T96}; see (\ref{ineq65}) for details. Note that  our inequality is a Gaussian bound. To the best of our knowledge, such type bounds for the maxima of empirical processes have not yet been obtained before.   The comparisons among our inequality with the known  inequalities in literature show that our concentration inequalities have certain advantage for small and moderate deviations $x$.

For the methodology, our method is based on Doob's martingale decomposition, Taylor's expansion and  Rio's inequality. It is similar to the method developed by Rio \cite{R09} for establishing deviation inequality for martingales. However, Rio's deviation inequality \cite{R09} does not provided  a Gaussian bound which plays an important role in this paper.

The paper is organized as follows. In Section \ref{sec02}, we present a Bernstein type condition for martingales and give a Gaussian version of  Freedman's  inequality under the stated condition.
 The applications are discussed in Sections \ref{sec3}, \ref{sec5} and \ref{sec7}.

\section{Concentration inequalities for martingales}\label{sec02}
Assume that we are given a sequence of real-valued martingale differences $%
(X _i,\mathcal{F}_i)_{i=0,...,n}, $ defined on some probability space $%
(\Omega ,\mathcal{F},\mathbf{P})$, where $X _0=0 $ and $\{\emptyset, \Omega\}=%
\mathcal{F}_0\subseteq ...\subseteq \mathcal{F}_n\subseteq \mathcal{F}$ are
increasing $\sigma$-fields. So, by definition, we have $\mathbf{E}[X_{i}|\mathcal{F}_{i-1}]= 0, \  i\in [1,n]$.
For any $j \in [1, n],$ denote by
$$
S_{0}=0,\ \ \ \ \ S_n=\sum_{i=1}^n X _i  \ \  \ \ \textrm{and} \ \ \ \  \langle S\rangle_n = \sum_{i=1}^n\mathbf{E}[X_{i}^2|\mathcal{F}_{i-1}].
$$
Hence $(S_i,\mathcal{F}_i)_{i=1,...,n}$ is a martingale.
Throughout the paper, $\textrm{Var}(S)$ stands for the variance of a random variable $S.$
It is well known that $\textrm{Var}(S_n)=\sum_{i=1}^n\mathbf{E}[ X_i^2 ].$
For simplicity of notation, denote   $$ \sigma^2=\frac{1}{n} \textrm{Var} (S_n).$$

\subsection{Gaussian version of  Freedman's  inequality}
In our main result, we shall make use of the following Bernstein type condition: there exists an  $\epsilon$, may depending on $n,$  such that
 \begin{eqnarray}\label{Btcondition}
\ \ \ \ \sum_{i=1}^n\mathbf{E}[|X_i|^k ] \leq \frac{1}{2} \frac{  k! \,\epsilon ^{k-2}  }{(k-1)^{k/2}}   \sum_{i=1}^n\mathbf{E}[X_i^2] \ \ \ \ \textrm{for all}\ k\geq  2.
\end{eqnarray}
 It is worth noting that the bounded martingale differences $|X_i| \leq a$ satisfy  condition (\ref{Btcondition}) with
 $\epsilon=2^{3/2}a/3$; see Theorem \ref{co1}. Moreover, the normal random variables also satisfy condition (\ref{Btcondition}); see (\ref{indu24}).
An equivalent condition of (\ref{Btcondition}) is given by the following theorem.
\begin{proposition}\label{th3s}
 Condition (\ref{Btcondition}) is equivalent to the following one: there exists a   $\rho$, such that
 \begin{eqnarray} \label{stBtdition}
\sum_{i=1}^n\mathbf{E}[(X_i)^k] \leq \frac{1}{2} \frac{  k! \, \rho^{k-2} }{(k-1)^{k/2}}  \sum_{i=1}^n\mathbf{E}[X_i^2],\ \ \ \ \textrm{for all even}\  k\geq 2.
\end{eqnarray}
\end{proposition}
\emph{Proof.}
 It is obvious that (\ref{Btcondition}) implies (\ref{stBtdition}) with $\rho=\epsilon.$  Next, we prove (\ref{stBtdition}) implies (\ref{Btcondition}). For even number $k= 2l, l\geq 1,$ it is obvious that  (\ref{stBtdition}) implies (\ref{Btcondition})   with $\epsilon=\rho.$  Thus we only need to prove the case of  $k= 2l + 1, l\geq 1.$   By   Cauchy-Schwarz's inequality, it is easy to see that
 \begin{eqnarray*}
\sum_{i=1}^n\mathbf{E}[|X_i|^k ] \leq \sum_{i=1}^n \Big(\mathbf{E}[|X_i|^{2l} ] \Big)^{1/2}  \Big(\mathbf{E}[|X_i|^{2(l+1)} ] \Big)^{1/2}.
\end{eqnarray*}
By  Cauchy-Schwarz's inequality again,  it follows that
 \begin{eqnarray*}
\sum_{i=1}^n \Big(\mathbf{E}[|X_i|^{2l} ] \Big)^{1/2}  \Big(\mathbf{E}[|X_i|^{2(l+1)} ] \Big)^{1/2} \leq  \Big(\sum_{i=1}^n \mathbf{E}[|X_i|^{2l} ] \Big)^{1/2}  \Big(\sum_{i=1}^n \mathbf{E}[|X_i|^{2(l+1)} ] \Big)^{1/2} .
\end{eqnarray*}
Thus
\begin{eqnarray*}
\sum_{i=1}^n\mathbf{E}[|X_i|^k ] \leq   \Big(\sum_{i=1}^n \mathbf{E}[|X_i|^{2l} ] \Big)^{1/2}  \Big(\sum_{i=1}^n \mathbf{E}[|X_i|^{2(l+1)} ] \Big)^{1/2} .
\end{eqnarray*}
Notice that $2l$ and $2(l+1)$ are even integers.
Using (\ref{stBtdition}), we obtain
 \begin{eqnarray*}
\sum_{i=1}^n\mathbf{E}[|X_i|^k] &\leq& \frac{1}{2}\ \rho^{2l-1} \bigg( \frac{(2l)! }{(2l-1)^l } \frac{(2l+2)! }{(2l+1)^{l+1} } \bigg)^{1/2} \sum_{i=1}^n\mathbf{E}[X_i^2]  \\
& =& \frac{1}{2} \frac{  k! \, \rho^{k-2} f(l) }{(k-1)^{k/2}}  \sum_{i=1}^n\mathbf{E}[X_i^2],
\end{eqnarray*}
where $$f(l)=\bigg(   \frac{  (2l)^{2l+1}(2l+2) }{ (2l-1)^l(2l+1)^{l+2} } \bigg)^{1/2}. $$
Since $f(l) \leq f(1)=(32/27)^{1/2}$ for all $l \geq 1,$ we have
 \begin{eqnarray*}
\sum_{i=1}^n\mathbf{E}[|X_i|^k] &\leq&  \frac{1}{2} \frac{  k! \, \rho^{k-2} f(1) }{(k-1)^{k/2}}  \sum_{i=1}^n\mathbf{E}[X_i^2]
\ \leq \ \frac{1}{2} \frac{  k! \, (f(1) \rho )^{k-2}  }{(k-1)^{k/2}}  \sum_{i=1}^n\mathbf{E}[X_i^2]
\end{eqnarray*}
for $k=2l+1, l\geq 1.$
Thus   (\ref{stBtdition}) implies (\ref{Btcondition})   with $\epsilon=(32/27)^{1/2}\rho.$\hfill\qed

In the following theorem, we give a Gaussian version of   Freedman's  inequality. Our inequality is similar to    Bernstein's  inequality.
\begin{theorem}
\label{th1} Assume condition (\ref{Btcondition}).  Then,
for any $0\leq t < \epsilon^{-1},$
 \begin{eqnarray}\label{sxds}
\mathbf{E}\bigg[\exp\bigg\{t \frac{ S_n}{\sqrt{n}} \bigg \}\bigg] \, \leq \,   \exp \Bigg\{ \frac{t^2 \sigma^2  }{2\, (1- t \epsilon)}    \Bigg\} ,
\end{eqnarray}
and,  for any positive $ x$,
\begin{eqnarray}
 \ \ \ \ \ \ \ \ \mathbf{P}\left( \max_{1\leq k\leq n} S_k \geq x \sqrt{\emph{Var}(S_n)}  \right) &\leq& \exp\bigg\{ - \frac{x^2}{  1+ \sqrt{ 1+2\, x \epsilon
  /\sigma   }  + x   \epsilon /\sigma   } \bigg\} \label{jnsk01}\\
  &\leq& \exp\bigg\{ - \frac{x^2}{  2 \big( 1+ x   \epsilon / \sigma   \big)} \bigg\}.\label{jnsk}
\end{eqnarray}
Moreover, the same inequalities hold when replacing $S_k$ by $-S_k.$
\end{theorem}

\begin{remark}
 Let us make some comments on Theorem  \ref{th1}.
\begin{enumerate}
\item Compared with  Freedman's inequality (\ref{freedmanxf2}),   our inequalities (\ref{jnsk01}) and (\ref{jnsk}) are expressed in terms of $\sigma^2$ instead of  $||\left\langle S\right\rangle_n||_\infty/n $. Moreover, our inequalities are valid for the martingales with unbounded differences.

\item Burkholder \cite{B73} proved that, for any $p>1,$
\begin{eqnarray}\label{f1.1s}
||S_n||_p \leq C_p || (X_1^2 + X_2^2+...+X_n^2 )^{1/2} ||_p .
\end{eqnarray}
In his paper in Ast\'{e}risque, Burkholder \cite{B88} obtained (\ref{f1.1s}) with $C_p=p-1$ for $p>2.$ He also proved that this constant
is optimal. From (\ref{f1.1s}), for any $p>2$ the following Marcinkiewicz-Zygmund type inequality holds:
\begin{eqnarray}\label{f1.2p}
||S_n||_p^2 \leq c_p (||X_1||_p^2 + ||X_2||_p^2+...+||X_n||_p^2 )
\end{eqnarray}
with $c_p=(p-1)^2.$   Rio \cite{R09} recently proved that (\ref{f1.2p}) holds with  $c_p= p-1,$  and that this constant cannot be  improved.

\item  Under the condition $\mathbf{E}[\exp\{X_k^2\}] \leq 1+ C$ for any positive $k$ and some positive  $C,$ Rio (cf.\ Corollary 2.2 of \cite{R09})  has obtained
the following result: for any positive $x,$
\begin{eqnarray}\label{fsfsf}
 \ \ \ \ \ \ \ \mathbf{P}\left( \max_{1\leq k\leq n} |S_k| \geq x \sqrt{ne/2} \right)  \leq  \frac{  C  }{ (2e)^{1/2}}   \ \big(  \cosh(x)-1 \big)^{-1} .
\end{eqnarray}
This bound has an exponentially decaying rate similar to that of (\ref{jnsk01}) as $x\rightarrow \infty$, but it
  is not a Gaussian bound.

\item  In certain cases, our inequality (\ref{fsf}) is  better than   Rio's  inequality (\ref{fsfsf}). For instance,
assume that  $|X_k|\leq \sqrt{k}$ for all $k\geq1,$ and that $  \emph{Var} (X_k), k\geq 1,$ are uniformly bounded. Thus
 $\mathbf{E}[\exp\{X_k^2\}] \leq   e^n$ for any positive $k$.
   Rio's  inequality (\ref{fsfsf}) implies that, for any positive $x,$
\begin{eqnarray}\label{cth}
\ \ \ \ \ \ \ \  \mathbf{P}\left( \max_{1\leq k\leq n} |S_k| \geq x \sqrt{ne/2} \right)  \leq  \frac{ e^{n}-1  }{ (2e)^{1/2}}    \ \big(  \cosh(x)-1 \big)^{-1} .
\end{eqnarray}
Our result implies that, for any positive $x,$
\begin{eqnarray*}
  \mathbf{P}\Big(\,  \max_{1\leq k\leq n} |S_k|\geq x \sqrt{n}  \Big)
   \leq  2\exp\left\{ - \frac{x^2}{2 (  \sigma^2 +    2^{3/2}  x /3 )} \right\}
\end{eqnarray*}
which is significantly better than (\ref{cth}).
\end{enumerate}
\end{remark}
\emph{Proof of Theorem \ref{th1}.}
 By Taylor's expansion of $e^x$ and $\mathbf{E}[S_n]=0$, we have, for all $t\geq 0,$
 \begin{eqnarray}\label{sxsa}
 \mathbf{E}\bigg[\exp\bigg\{t \frac{ S_n}{\sqrt{n}} \bigg \}\bigg] = 1+ \sum_{k=2}^{\infty} \frac{t^k}{k!}  \mathbf{E}\Big[\Big(\frac{ S_n}{\sqrt{n}}\Big)^k \Big].
\end{eqnarray}
Using the following Rio inequality (see Theorem 2.1 of \cite{R09}): for any $p\geq 2,$
\begin{eqnarray}\label{R}
\Big(\mathbf{E}[|S_n|^p] \Big)^{2/p} \leq  (p-1)   \Big(\sum_{i=1}^n  \big(\mathbf{E}[|X_i|^p] \big)^{2/p} \Big),
\end{eqnarray}
we get, for all $k\geq 2,$
\begin{eqnarray}\label{ssxsfsd}
\mathbf{E}[|S_n|^k] \leq (k-1)^{k/2} \Big(\sum_{i=1}^n \big(\mathbf{E}[|X_i|^k] \big)^{2/k}\Big)^{k/2}.
\end{eqnarray}
Hence, by the inequality
\begin{eqnarray}\label{fnlbjs}
 (a_1+a_2+...,+a_n)^p \leq n^{p-1}(a_1^p+a_2^p+...,+a_n^p),
\end{eqnarray}  inequality (\ref{ssxsfsd}) implies that, for all $k\geq 2,$
\begin{eqnarray}\label{sddd}
\mathbf{E}[|S_n|^k] \leq (k-1)^{k/2} n^{k/2-1} \sum_{i=1}^n  \mathbf{E}[|X_i|^k]  .
\end{eqnarray}
Applying the last inequality to (\ref{sxsa}), we obtain
 \begin{eqnarray}\label{bnlts}
\mathbf{E}\bigg[\exp\bigg\{t \frac{ S_n}{\sqrt{n}} \bigg \}\bigg] \leq 1+ \sum_{k=2}^{\infty} \Big( \frac{t^k}{k!}  (k-1)^{k/2} n^{ -1 } \sum_{i=1}^n  \mathbf{E}[|X_i|^k] \Big).
\end{eqnarray}
Hence  condition (\ref{Btcondition}) implies that, for all $0\leq t < \epsilon^{-1}$,
 \begin{eqnarray}\label{ineq23}
\mathbf{E}\bigg[\exp\bigg\{t \frac{ S_n}{\sqrt{n}} \bigg \}\bigg] \leq 1+ \sum_{k=2}^{\infty}  \frac{ \ \sigma^2}{2 } \, t^k \epsilon^{k-2} = 1+   \frac{t^2 \sigma^2  }{2\,(1- t \epsilon)} .
\end{eqnarray}
By the inequality $1+x \leq e^x  ,$ it follows that, for all $0\leq t < \epsilon^{-1}$,
 \begin{eqnarray*}
\mathbf{E}\bigg[\exp\bigg\{t \frac{ S_n}{\sqrt{n}} \bigg \}\bigg] \ \leq \   \exp \Bigg\{ \frac{t^2  \sigma^2  }{2\,(1- t \epsilon)}   \Bigg\} ,
\end{eqnarray*}
which gives the desired inequality (\ref{sxds}).  Applying  Doob's maximal inequality   to the nonnegative submartingale $(\exp \{t S_k /(\sqrt{n} \sigma)\}, \mathcal{F}_k)_{1\leq k \leq n}, 0\leq t < \epsilon^{-1},$  it is easy to see that, for any $ x \geq 0,$
\begin{eqnarray*}
  \mathbf{P}\left( \max_{1\leq k\leq n} S_k \geq x \sqrt{n} \sigma \right) &\leq&  \mathbf{P}\left( \max_{1\leq k\leq n} \exp\bigg\{t \frac{ S_k}{\sqrt{n}\sigma} \bigg \}   \geq  \exp\Big\{t x   \Big \}  \right) \\
  &\leq&  \exp\Big\{ - t x \Big \} \mathbf{E}\left[   \exp\bigg\{t \frac{ S_n}{\sqrt{n} \sigma } \bigg \} \right].
\end{eqnarray*}
Hence
\begin{eqnarray*}
  \mathbf{P}\left( \max_{1\leq k\leq n} S_k \geq x \sqrt{n} \sigma \right)
    &\leq& \inf_{0\leq t < \epsilon^{-1}}  \exp\bigg\{ - t x  + \frac{t^2    }{2\,(1- t \, \epsilon/ \sigma)}    \bigg \} \\
    &=& \exp\bigg\{ - \frac{x^2}{  1+ \sqrt{ 1+2\, x \epsilon
  /\sigma   }  + x   \epsilon /\sigma   } \bigg\},
\end{eqnarray*}
which gives   (\ref{jnsk01}).   Using the   inequality
$   \sqrt{ 1+2\, x \epsilon / \sigma    }  \leq  1+ \, x \epsilon / \sigma ,$
we get  (\ref{jnsk})  from (\ref{jnsk01}).  \qed

 From (\ref{jnsk}), it is interesting to  see that for small $x$, i.e.\ $0\leq x =o(\sigma /\epsilon),$ the tail probabilities on partial maximum martingales $\mathbf{P}\big( \max_{1\leq k\leq n} S_k \geq x \sqrt{\textrm{Var}(S_n)} \ \big)$ have
the exponential Gaussian bound $\exp\big\{- x^2/2   \big\}.$
We wonder the range $0\leq x =o(\sigma /\epsilon)$ can be extended to a larger one under the stated condition. In other words,  could we establish the following inequality under  condition (\ref{Btcondition}). There  exist   some positive constants $b_n,$   with $b_n \rightarrow 0$ as $n\rightarrow\infty,$ such that,  for any positive $ x$,
\begin{eqnarray}\label{thn3s}
  \ \ \mathbf{P}\Big( \max_{1\leq k\leq n} S_k \geq x \sqrt{\textrm{Var}(S_n)}  \Big)
  \leq \exp\bigg\{ - \frac{x^2}{  2 \big( 1 + x \,  b_n \epsilon/\sigma    \big)} \bigg\}\  ?
\end{eqnarray}
When $(X_i)_{i\geq1}$ are i.i.d.,  Bernstein's inequality shows that we can take  $b_n=1/\sqrt{n}.$ However, for martingale differences,
the following theorem  gives a negative answer to question (\ref{thn3s}).
\begin{theorem}\label{th31}
There exists a stationary sequence of martingale differences  $(X_{i}, \mathcal{F}_{i})_{i\geq 1}$ satisfying condition (\ref{Btcondition})
and, for all $x  \rightarrow \infty$,
\begin{eqnarray}\label{fnk}
\mathbf{P}\Big(  \max_{1\leq k \leq n}  S_k  \geq  x\sqrt{\emph{Var}(S_n)} \Big)  \geq    \exp \Big\{ -C\,  x \Big\},\  \ \
\end{eqnarray}
where $C$ is a positive constant and does not depend on $n$ and $x$.
\end{theorem}

Assume that there exist two positive constants $a, b$ such that  $ \sigma /\epsilon \in [a, b]$ uniformly for all $n\geq 1$, which holds for stationary martingale difference sequence. The inequality (\ref{thn3s}) implies that, for $x_n \rightarrow \infty$ and $x_n=o(b_n^{-1} )$ as $n\rightarrow \infty,$
\begin{eqnarray}\label{g63s}
  \mathbf{P}\Big( \max_{1\leq k\leq n} S_k \geq x_n \sqrt{\textrm{Var}(S_n)} \Big)
  \leq \exp\Big\{ - C\,  x_n^2  \Big\},
\end{eqnarray}
where $C$ is a positive constant and does not depend on $n.$
For $x_n \rightarrow \infty$ and $x_n=o(b_n^{-1})$ as $n\rightarrow \infty,$ it is obvious that $\exp\big\{ - C_1 \, x_n^2  \big\}  < \big\{ -C_2\,  x_n \big\},$ and thus (\ref{g63s}) contradicts (\ref{fnk}). Hence, we cannot establish  (\ref{thn3s}) under  condition (\ref{Btcondition}) even for stationary martingale difference sequences.

\noindent\emph{Proof of Theorem \ref{th31}}.
  We proceed as in Fan \textit{et al.}\ \cite{FGL15}. Let $X$  be a standard normal random variable.
By Stirling's formular
$$
n!= \sqrt{2\pi n}\, n^ne^{-n}e^{\frac{1}{12 n\theta_n} }, \ \ \ \ 0 \leq \theta_n \leq 1,
$$
it is easy to verify that, for all $k \geq3,$
 \begin{eqnarray}
  \mathbf{E}[|X| ^k ] &\leq& (k-1)!! \leq \sqrt{k!} \nonumber \\
 &\leq &\frac{1}{2}\sqrt{k! \,}  \Big(1+ \frac{1}{k-1}\Big)^{\frac k 2}  e^{\frac{1}{24k \theta_k} -\frac k2}\Big(\frac{2}{e^{1/2}}\Big)^{k-2} \nonumber \\
 &  \leq&  \frac{1}{2} \frac{  k! \, }{(k-1)^{k/2}}  \Big(\frac{2}{e^{1/2}}\Big)^{k-2}.\label{indu24}
\end{eqnarray}
 Assume that $(\xi_i)_{i\geq 1}$ are Rademacher random variables independent of $X$, i.e. $\mathbf{P}(\xi_i=1)=\mathbf{P}(\xi_i=-1)=\frac12$. Set $X_i=X\xi_i$ and $\mathcal{F}_i=\sigma(X, (\xi_k)_{k=1,...,i})$. Then by Theorem  \ref{th3s},   $(X_i, \mathcal{F}_i)_{i\geq 1}$ is a stationary sequence of martingale differences and satisfies (\ref{Btcondition}) with $\sigma=1$ and $\epsilon=2/\sqrt{e}.$  It is easy to see that, for any positive $x,$
\begin{eqnarray*}
\mathbf{P}\Big(  \max_{1\leq k \leq n}  S_k \geq x \sqrt{n}  \Big)   \geq   \mathbf{P}\Big(    S_n \geq x \sqrt{n}  \Big)
  \geq    \mathbf{P}\Big(  \sum_{i=1}^n \xi_i \geq \sqrt{n x } \Big)  \mathbf{P} \Big(   X \geq \sqrt{x}   \Big).
\end{eqnarray*}
Since  for $x$ large enough,
\[
\mathbf{P}\Big(  \sum_{i=1}^n \xi_i \geq  \sqrt{n x} \Big) \geq   \exp \left\{ - \frac{ (\sqrt{n x})^2}n   \right\}=e^{- x},
\]
we get
\begin{eqnarray}
\mathbf{P}\Big(  \max_{1\leq k \leq n}  S_k  \geq  x\sqrt{n} \Big)  \geq  \frac{1}{ \sqrt{2\pi}(1+ \sqrt{x})   } \exp\Big\{-x  - \frac12 x \Big\}
 \geq    \exp \Big\{ -2\, x \Big\}, \nonumber
\end{eqnarray}
which gives (\ref{fnk}). This ends the proof of Theorem \ref{th31}.\hfill\qed

\subsection{Case for bounded differences} \label{sec2}
When the martingale differences are  bounded, Theorem  \ref{th1} implies the following corollary which is even better than the classical Azuma-Hoeffding inequality for moderate deviations $x$.
\begin{theorem}
\label{co1} Assume that  $|X_i|\leq a $ for some constant $a$ and all $i \in [1, n]$.
Then the inequalities (\ref{sxds}), (\ref{jnsk01}) and (\ref{jnsk}) hold with  $\epsilon= 2^{3/2}a/3.$
Moreover, the same inequalities hold   when replacing $S_k$ by $-S_k.$
\end{theorem}
\emph{Proof}.   Define the function $$f(k)=\frac{1}{2} k! \frac{1}{(k-1)^{k/2}}$$ for $k \geq 3$.
Note that the function $g(k)=f(k)(f(3))^{2-k}$ is increasing in $k\geq 3$ and $g(3)=1$. Thus
$f(3)^{k-2}\leq f(k) $ for any $k\geq 3$.
If $|X_i|\leq a$,
 then, for any $k\geq 3,$
\begin{eqnarray}
\sum_{i=1}^n\mathbf{E}[|X_i|^k]  &\leq &  (f(3))^{k-2} \Big(\frac{a}{f(3)} \Big)^{k-2}    \sum_{i=1}^n\mathbf{E}[X_i^2] \nonumber \\
&\leq&  f(k)  \Big(\frac{a}{f(3)} \Big)^{k-2}     \sum_{i=1}^n\mathbf{E}[X_i^2].
\end{eqnarray}
Hence condition (\ref{Btcondition}) holds with $\epsilon =a/ f(3) = 2^{3/2}a /3 $.
By Theorem \ref{th1}, we obtain the desired inequalities. \hfill\qed

Next, we examine the general $\mathcal{F}_n$-measurable random  variables.
\begin{theorem}
\label{co2}
Let $S_n$ be a random variable, and let $\{\emptyset, \Omega\}=\mathcal{F}_0\subset \mathcal{F}_1 \subset ... \subset \mathcal{F}_n$  be a sequence of increasing $\sigma$-fields. Assume that
$S_n$ is $\mathcal{F}_n$-measurable. Assume that there exists a sequence of $\mathcal{F}_{i}$-measurable random variables $T_i$ such  that
\begin{eqnarray}\label{fskklvcn}
T_{i-1} \leq  \mathbf{E}[S_n | \mathcal{F}_{i}] \leq  T_{i-1} + M , \ \ \ \ i=1,...,n,
\end{eqnarray}
where $M$ is a constant.
 Then,   for any positive $ x$,
\begin{eqnarray}\label{ffg}
  \ \ \ \ \ \ \ \ \ \mathbf{P}\Big(\, S_n-\mathbf{E}[S_n]\geq x  \Big)
   \leq  \exp\left\{ - \frac{x^2}{2 ( \emph{Var}(S_n) +    2^{3/2}M \sqrt{n} x /3 )} \right\}.
\end{eqnarray}
Moreover, the same inequality holds when replacing $S_n$ by $-S_n.$
\end{theorem}


Comparing with   Azuma-Hoeffding's inequality,  our inequality (\ref{ffg}) has certain advantages.
Assume  that
\begin{eqnarray}\label{fskklvfcn}
T_{i-1} \leq  \mathbf{E}[S_n | \mathcal{F}_{i}] \leq  T_{i-1} + M_i , \ \ \ \ i=1,...,n,
\end{eqnarray}
where $T_i$ are $\mathcal{F}_{i}$-measurable random variables and $M_i$ are some constants. Under condition (\ref{fskklvfcn}), it is well known that $\textrm{Var}(S_n) \leq \frac14\sum_{i=1}^n M_i^2; $ see McDiarmid \cite{M89}.
 The classical Azuma-Hoeffding inequality states that,
for any positive $x,$
\begin{eqnarray}\label{newamb}
 \mathbf{P}\Big(S_n-\mathbf{E}[S_n]  \geq x  \Big)
 \leq  \exp \left\{-\frac{ 2 \, x^2}{ \sum_{i=1}^n M_i^2} \right\} .
\end{eqnarray}
First, it is easy to see that inequality (\ref{ffg}) is a Gaussian  bound, while   Azuma-Hoeffding's inequality
does not share this feature.
Second,  inequality (\ref{ffg}) is better than  Azuma-Hoeffding's inequality for all $x$ in the range
$$0 \leq x \leq \frac{3}{2^{3/2}M \sqrt{n}}\Big( \frac14\sum_{i=1}^n M_i^2 - \textrm{Var}(S_n)  \Big) ,$$
where $M= \max\{M_i: i=1,...,n\}.$
To illustrate the last range, consider the case that $ \textrm{Var} (S_n)$ is in order of $n$ and $M_i= i^\alpha$ for some $ \alpha> 0$ and all  $ i\geq 1.$   Our inequality (\ref{ffg})  with $M=n^\alpha$  improves  Azuma-Hoeffding's inequality (\ref{newamb}) for all $$0\leq x \leq \frac{3}{2^{3/2} n^{\alpha + 1/2}  }\Big(   \frac14 \sum_{i=1}^n i^{2\alpha}  - \textrm{Var}(S_n) \Big).$$
Notice that for the right hand side of the last inequalities, it holds   $$ \frac{1}{ n^{\alpha + 1/2}  }\Big( \frac14  \sum_{i=1}^n i^{2\alpha}  - \textrm{Var}(S_n) \Big)\rightarrow \frac{1}{8\alpha+4}n^{\alpha+ \frac12} , \ \ \ \ n \rightarrow \infty.$$
Recall that $ \textrm{Var} (S_n) $ is in order of $n.$ Thus inequality (\ref{fsf}) improves   Azuma-Hoeffding's inequality (\ref{ineq3}) for all \emph{standard} $x$ in a range $0\leq x =O(n^\alpha), n\rightarrow \infty. $ This range is quite large.
Here and after, call $x$ is  standard   if we refer to the  tail probabilities for the standardized sums, i.e.\ $\mathbf{P}\big( S_n /\sqrt{\textrm{Var} (S_n)} \geq x \big).$


\noindent\emph{Proof of Theorem \ref{co2}}.
Let $  S_n =\sum_{i=1}^n X_i$ be Doob's martingale decomposition of $S_n,$
where
\begin{eqnarray*}
X_i = \mathbf{E}[S_n | \mathcal{F}_{i}]-\mathbf{E}[S_n | \mathcal{F}_{i-1}].
\end{eqnarray*}
Thus (\ref{fskklvcn}) implies that
\begin{eqnarray}\label{finew32}
T_{i-1}   -\mathbf{E}[S_n | \mathcal{F}_{i-1}] \leq  X_i \leq T_{i-1} + M  -\mathbf{E}[S_n | \mathcal{F}_{i-1}].
\end{eqnarray}
By the fact $\mathbf{E}[X_i | \mathcal{F}_{i-1}]=0,$ it follows that
   $$T_{i-1}   -\mathbf{E}[S_n | \mathcal{F}_{i-1}] \leq 0\ \ \ \ \textrm{and}\ \ \ \  T_{i-1} + M  -\mathbf{E}[S_n | \mathcal{F}_{i-1}] \geq 0.$$
Thus we have
$$T_{i-1} + M  -\mathbf{E}[S_n | \mathcal{F}_{i-1}] \leq M  \ \ \ \ \textrm{and}\ \ \ \   -M \leq T_{i-1}   -\mathbf{E}[S_n | \mathcal{F}_{i-1}].$$
Returning to (\ref{finew32}), we get $|X_i| \leq M.$ Then Theorem \ref{co2} follows from Theorem \ref{co1}.\hfill\qed

\section{Concentration for functions of independent random variables}\label{sec3}
In this section we give  some   applications of our inequalities to functions of independent random variables.
Let $(\xi_i)_{i=1,...,n}$ be a sequence of independent random variables with values
in some complete separable metric space $\mathcal{X}$.  Let $f$ be a function from $\mathcal{X}^n$ to $\mathbf{R}.$
We are interested in the concentration for the random function $f(\xi_1,...,\xi_n)$.

\subsection{McDiarmid type inequality}
Assume that $(\mathcal{X}, d_i), i=1,...,n,$ are separable metric
spaces with  positive finite diameters $M_i,$ i.e.\ $d_i(\cdot, \cdot)\leq M_i,$  where $d_i$ are some distances on $\mathcal{X}.$ Let $f$ be a separately Lipschitz function  such that
\begin{eqnarray}\label{spcond}
\ \ \ \ |f(x_1, . . . , x_n)-f(y_1, . . . , y_n)| \, \leq \,   d_1(x_1, y_1)+ ...  +  d_n(x_n, y_n).
\end{eqnarray}
Set
\begin{eqnarray}\label{mdicss}
Z_n= f( \xi_1,..., \xi_n ).
\end{eqnarray}
 McDiarmid  \cite{M89}  has obtained  the following concentration inequality: for any positive  $t$,
\begin{eqnarray}\label{fsfMCDi}
\mathbf{P}\Big( Z_n - \mathbf{E}[Z_n] \geq t\sqrt{n} \Big) \ \leq \ \exp\Big\{-\frac{2 n t^2}{  T_n^2} \Big\},
\end{eqnarray}
where
$$ T_n^2=\sum_{i=1}^n M_i^2.$$
He also showed that  $\textrm{Var}  (Z_n) \leq \frac14 T_n^2 .$
 \mbox{McDiarmid's} inequality can be regarded as a generalization of   Azuma-Hoeffding's inequality in the functional setting.
 Recently, Rio \cite{R13} has obtained the following improvement on  \mbox{McDiarmid's} inequality: for any $t \in [0, 1],$
\begin{eqnarray}\label{Rios}
\mathbf{P}\Big( Z_n - \mathbf{E}[Z_n] \geq t D_n  \Big) \ \leq \ (1-t)^{t(2-t)D_n^2/T_n^2 },
\end{eqnarray}
where $$D_n=\sum_{i=1}^n M_i. $$
As pointed out by Rio, his bound (\ref{Rios}) is less than \mbox{McDiarmid's} bound (\ref{fsfMCDi}).
Moreover, Rio's bound has an interesting feature, that is
$ \mathbf{P}\big( Z_n - \mathbf{E}[Z_n] \geq t \big)=0$ when $t> D_n$, which coincides with the property $Z_n - \mathbf{E}[Z_n] \leq D_n.$
Rio's inequality (\ref{Rios}) can be rewritten in the following form: for any $t \in [0, D_n/\sqrt{n}\,],$
\begin{eqnarray}\label{Riods}
\ \   \ \ \ \ \ \mathbf{P}\Big( Z_n - \mathbf{E}[Z_n] \geq t  \sqrt{n} \Big)  \leq \exp\bigg\{ t \sqrt{n} \Big(2-t\frac{\sqrt{n}}{D_n}\Big)\frac{ D_n }{T_n^2} \ln\Big(1-t\frac{ \sqrt{n}}{D_n} \Big)  \bigg\} .
\end{eqnarray}
Notice that for all $0\leq t =o(D_n/\sqrt{n}),$ Rio's bound (\ref{Riods}) is of type $\exp\big\{-\frac{2 n t^2}{  T_n^2} (1+o(1))\big\},$ which is similar to McDiarmid's bound (\ref{fsfMCDi}).

In the following theorem, we extend   Theorem \ref{co1} to  the functional setting.
\begin{theorem} \label{co3}
Let $Z_n$ be defined by (\ref{mdicss}). Denote  $$\sigma^2=\frac{1}{n}  \emph{Var}(Z_n) \  \ \ \ \ \textrm{and }\ \ \ \ \ M=\max\{M_i:\,  1\leq i \leq n \} .$$ Then,  for any positive $ t$,
\begin{eqnarray}
  \mathbf{P}\Big(\, Z_n-\mathbf{E}[Z_n]\geq t \sqrt{\emph{Var}(Z_n)} \Big)  &\leq &   \exp\Bigg\{ - \frac{t^2}{  1+ \sqrt{1+\frac{2^{5/2}  }{3 } M t /\sigma } \, + \frac{ 2^{3/2} }{ 3  } M t/\sigma  } \Bigg\} \label{McdiamidP} \\
  &\leq&  \exp\bigg\{ - \frac{t^2}{2 \big( 1 +  \frac{2^{3/2}}{3} M t /\sigma  \big)} \bigg\}.\label{McdiamidR}
\end{eqnarray}
Moreover, the same inequalities hold when replacing $Z_n$ by $-Z_n.$
\end{theorem}
\emph{Proof.} Denote $\mathcal{F}_i = \sigma\{ \xi_j,\ 1\leq j \leq i \}. $  Let $Z_n-\mathbf{E}[Z_n] = S_n$ be Doob's martingale decomposition of $Z_n,$ where
$
X_i = \mathbf{E}[Z_n | \mathcal{F}_{i}]-\mathbf{E}[Z_n | \mathcal{F}_{i-1}].
$
Let $(\xi_i')_{i=1,..,n}$ be an independent copy of the random variables $(\xi_i)_{i=1,..,n}.$
Then it is easy to see that
\begin{eqnarray}
X_i &=& \mathbf{E}[ f( \xi_1,... , \xi_i, \xi_{i+1}' , ..., \xi_n' ) - f( \xi_1,...,  \xi_{i-1}, \xi_{i}' ,..., \xi_n ')   | \mathcal{F}_{i}]  \nonumber\\
&\leq& \mathbf{E}[ d_i( \xi_i, \xi_i' ) | \mathcal{F}_{i} ]  \label{ineq27}\\
&\leq& M. \nonumber
\end{eqnarray}
Similarly, we have $ X_i  \geq - M.$ Thus $|X_i|\leq M.$  Then the inequalities (\ref{McdiamidP}) and (\ref{McdiamidR}) follow by Theorem
\ref{co1}. This completes the proof of theorem. \hfill\qed

Comparing to  the inequalities of  McDiarmid  (\ref{fsfMCDi}) and Rio (\ref{Riods}), our inequalities (\ref{McdiamidP}) and (\ref{McdiamidR}) have
the following two interesting features. The first feature is that for $0\leq t =o(1),$ our bounds have the exponential Gaussian form $\exp\big\{- t^2/ 2 \big\},$ in contrast to the bounds of McDiarmid and Rio which do not share this property.
The second feature is that our inequalities are better than  (\ref{fsfMCDi}) and  (\ref{Riods})  in the range $0\leq t = \frac{3}{ 2^{3/2}M }\big(\frac1{4n}T_n^2 -\sigma^2 \big)$.
 Recall that when  $M_i=i^\alpha, \alpha>0,$ for all $i$ and $\sigma$ is uniformly bounded for $n$, then the last range reduces to $0\leq t     =O(n^\alpha),   n\rightarrow \infty.$
Notice that for all $0\leq t =o(D_n/\sqrt{n}),$ Rio's bound (\ref{Riods}) is similar to McDiarmid's bound.

Notice that (\ref{McdiamidR}) is better than McDiarmid's inequality for moderate deviations $t$, and that \mbox{McDiarmid's} inequality  is better than (\ref{McdiamidR}) for $t$ large enough. Therefore, it is natural to minimize the two bounds. Then we derive  the following corollary.
\begin{corollary}\label{co3s} Assume  conditions of Theorem   \ref{co3}.  Then,  for any positive $t,$
\begin{eqnarray}\label{fdv}
  \ \ \ \ \ \ \mathbf{P}\Big(\, Z_n-\mathbf{E}[Z_n]\geq t \sqrt{n}  \Big)
   \leq  \exp\bigg\{ - \frac{t^2}{2 \min\big\{  \sigma^2 +  \frac{2^{3/2}}{3} M t ,  \frac1{4 n}  T_n^2 \big\}} \bigg\}.
\end{eqnarray}
Moreover, the same inequality holds when replacing $Z_n$ by $-Z_n.$
\end{corollary}

\begin{remark}
For the sake of simplicity  we restrict ourselves to the slightly weaker bound (\ref{fdv}), although (\ref{fdv}) can be improved
by minimizing Rio's bound  (\ref{Rios}) and the bound (\ref{McdiamidP}).
\end{remark}

In many applications, we would like to obtain information about the variance $\textrm{Var}  (Z_n).$   Thus we collect the following three estimations of $\textrm{Var}  (Z_n)$, where the last two estimations of $\textrm{Var}  (Z_n)$ can be found in Boucheron, \mbox{Lugosi} and Bousquet \cite{BLB}.
\begin{enumerate}

\item Let $(\xi_i')_{i=1,..,n}$ be an independent copy of the random variables $(\xi_i)_{i=1,..,n}.$ Write
  $$ Z_i' = f( \xi_1,...,\xi_{i-1} ,\xi_i',\xi_{i+1},..., \xi_n ).$$ The Efron-Stein inequality (cf.\ Efron-Stein \cite{es} and Steele \cite{st}) states that
\begin{eqnarray}
 \textrm{Var}  (Z_n)  \leq \frac12 \sum_{i=1}^n \mathbf{E}[ (Z_n-Z_i' )^2  ].
\end{eqnarray}
  In particular, under condition (\ref{spcond}), it  implies that $\textrm{Var} (Z_n) \ \leq \ \frac12  \sum_{i=1}^n\mathbf{E}[ d_i^2( \xi_i, \xi_i' )   ].$

\item Write $\mathbf{E}_i$ for the expected value without respect to the variable $\xi_i,$ that is $$\mathbf{E}_i = \mathbf{E}[ \,\cdot\, |\xi_1,...,\xi_{i-1},\xi_{i+1},...,  \xi_n ].$$ Then
      $$ \textrm{Var}  (Z_n)  \leq \sum_{i=1}^n \mathbf{E}[ (Z_n-\mathbf{E}_i [Z_n] )^2  ].$$

\item Write
  $$ \widetilde{Z}_i  = f_i(\xi_1,...,\xi_{i-1},\xi_{i+1},...,  \xi_n)$$
  for arbitrary measurable function $f_i: \mathcal{X}^{n-1} \rightarrow \mathbf{R}$ of $n-1$ variables.  Then
\begin{eqnarray}
\textrm{Var}  (Z_n)  \leq   \sum_{i=1}^n \mathbf{E}[ (Z_n-\widetilde{Z}_i  )^2  ]. \label{fdfgfp}
\end{eqnarray}
\end{enumerate}

\subsection{Concentration for self-bounding functions}
Let $\xi \in \mathcal{X}^n.$ Denote    $\xi^{(i)}=(\xi_1,...,\xi_{i-1},\xi_{i+1},..,\xi_n) \in \mathcal{X}^{n-1},$  obtained by dropping the $i$-th component of $\xi.$ For each $i\leq n,$ denote by $f_i$ a measurable function
 from $\mathcal{X}^{n-1}$ to $\mathbf{R}.$
Set
\begin{eqnarray}\label{fss57}
Z= f( \xi_1,..., \xi_n ).
\end{eqnarray}
We have the following concentration inequalities for  $Z$ around its expected value.
\begin{theorem}\label{th1ij} Define $Z$ by (\ref{fss57}). Assume that
\begin{eqnarray}\label{coind01}
0\leq f(\xi)-f_i(\xi^{(i)}) \leq 1
\end{eqnarray}
for all $i=1,...,n$ and
all $\xi \in \mathcal{X}^n.$ Denote   $\sigma^2=\frac1n \emph{Var}(Z).$
Then, for any positive $t$,
\begin{eqnarray}
 \mathbf{P}\Big(\, Z-\mathbf{E}[Z]\geq t \sqrt{ \emph{Var}(Z) } \ \Big)
   &\leq&  \exp\Bigg\{ - \frac{t^2}{  1+ \sqrt{1+\frac{2^{5/2}  }{3 }   t /\sigma  }  +  \frac{2^{3/2}}{3} t /\sigma    } \Bigg\} \label{ineq46}\\
  &\leq&  \exp\Bigg\{ - \frac{t^2}{2 \big(  1 +  \frac{2^{3/2}}{3}  t /\sigma \big)} \Bigg\}. \label{ineq47}
\end{eqnarray}
Moreover, the same inequalities hold  when replacing $Z$ by $-Z.$
\end{theorem}

A function $f$ is called $(a, b)$-self-bounding, introduced by McDiarmid and Reed \cite{MR0}, if condition (\ref{coind01}) holds and, moreover,  for some $a>0, b\geq0$ and
all $\xi \in \mathcal{X}^n,$
\begin{eqnarray} \label{conuml}
\sum_{i=1}^{n}\Big( f(\xi)-f_i(\xi^{(i)})\Big) \leq af(\xi)+b.
\end{eqnarray}
In particular,  $(1, 0)$-self-bounding function is known as self-bounding function; see Boucheron, Lugosi and Massart \cite{B09}.
For any $(a, b)$-self-bounding  function $f,$  McDiarmid and Reed  \cite{MR0} proved that, for any positive  $t,$
\begin{eqnarray}\label{ffs}
  \mathbf{P}\Big(\, Z-\mathbf{E}[Z]\geq t   \Big)
   \leq  \exp\bigg\{ - \frac{t^2}{2 \big( a \mathbf{E}[Z]  +b   +  c \, t   \big)} \bigg\},
\end{eqnarray}
where $c=a$ if $Z=f( \xi_1,..., \xi_n )$, and $c=1/3$ if $Z=-f( \xi_1,..., \xi_n ).$  See also Boucheron, Lugosi and Massart   \cite{B09}
  for the self-bounding functions.

  It is   worth noting that  (\ref{ineq47}) does not assume condition (\ref{conuml}). Hence, our inequality (\ref{ineq47}) extends
  the  inequality of McDiarmid and Reed \cite{MR0}.
  Moreover,   for
$(a, b)$-self-bounding  function, it holds  $ \textrm{Var}(Z) \leq a \mathbf{E}[Z]  +b .$
Indeed, by (\ref{fdfgfp}), (\ref{coind01}) and (\ref{conuml}), it is easy to see that
\begin{eqnarray}\label{inbeg59}
\textrm{Var}  (Z  )&\leq&  \mathbf{E}\Big[ \sum_{i=1}^n (Z -\widetilde{Z}_i  )^2  \Big]\ \leq \ \mathbf{E}\Big[ \sum_{i=1}^n (Z -\widetilde{Z}_i  )  \Big]  \ \leq \  \mathbf{E}[  aZ +b ]\nonumber \\
&=&a \mathbf{E}[Z]  +b ,
\end{eqnarray}
where $\widetilde{Z}_i =f_i( \xi^{(i)}).$
Thus, our bound (\ref{ineq47}) is less than the bound  of   McDiarmid and Reed (\ref{ffs}) for standard   $0\leq t =O\Big( \sqrt{ \frac{n}{\textrm{Var}  (Z  )} }\Big( \frac{a\mathbf{E}[Z]  +b}{\textrm{Var}  (Z  )}-1 \Big) \Big)$ as $n\rightarrow\infty$.
 The last range is large provided that $ n /  \textrm{Var}  (Z  )  \rightarrow \infty, n\rightarrow \infty$.

 \noindent\emph{Proof of Theorem \ref{th1ij}}.  Let $(\mathcal{F}_i)_{i=1,..,n}$ be the  natural filtration of the random variables $(\xi_i)_{i=1,..,n}.$  Denote  $Z=f(\xi),$ where $\xi=( \xi_1,..., \xi_n).$ Let $Z-\mathbf{E}[Z] = S_n$ be Doob's martingale decomposition of $Z,$ where $X_i = \mathbf{E}[Z | \mathcal{F}_{i}]-\mathbf{E}[Z | \mathcal{F}_{i-1}].$
By condition (\ref{coind01}), it is easy to see that
\begin{eqnarray}
X_i &\leq& \mathbf{E}[ 1+ f_i( \xi^{(i)}) |  \mathcal{F}_{i}]   -\mathbf{E}[  f( \xi)   | \mathcal{F}_{i-1}]  \nonumber\\
&=&  \mathbf{E}[ 1+ f_i( \xi^{(i)}) - f( \xi)     | \mathcal{F}_{i-1}] \nonumber\\
&\leq& 1. \nonumber
\end{eqnarray}
Similarly, we have
\begin{eqnarray}
X_i &\geq& \mathbf{E}[  f_i( \xi^{(i)}) |  \mathcal{F}_{i}]   -\mathbf{E}[  f( \xi)   | \mathcal{F}_{i-1}]  \nonumber\\
&=&  \mathbf{E}[   f_i( \xi^{(i)}) - f( \xi)     | \mathcal{F}_{i-1}] \nonumber\\
&\geq& -1. \nonumber
\end{eqnarray} Thus $|X_i|\leq 1.$ Then the inequalities (\ref{ineq46}) and (\ref{ineq47}) follow from Theorem
\ref{co1}.
When $Z=-f(\xi),$  the proof is similar to the case of $Z=f(\xi)$.    \hfill\qed

\subsection{Vapnik-Chervonenkis entropies}
Let $\mathcal{A}$ be an arbitrary collection of subsets of $\mathcal{X},$ and let
$x=(x_1,...,x_n)$ be a vector of $n$ points of $\mathcal{X}.$ Define the trace of $\mathcal{A}$
on $x$ by
$$\textrm{tr}(x)=\{  A\cap \{x_1,...,x_n\}: \ A \in \mathcal{A}\}.$$
The  \emph{shatter coefficient}  of $\mathcal{A}$ in $x$ is defined by $T(x)=|\textrm{tr}(x)|,$ namely the size of the trace.
$T(x)$ is the number of different subsets of the $n$-point set $\{x_1,...,x_n\}$ generated by
intersecting it with elements of $\mathcal{A}.$
The Vapnik-Chervonenkis  entropy (VC entropy) is defined as $$H(x)=\log_2T(x),$$
 with the convention that   $\log_20=-1.$
Note that
\begin{eqnarray}\label{coblm}
0\leq H(x)-H(x^{(i)})  \leq 1.
\end{eqnarray} The VC entropy $H(x)$ is of
particular interest, as it plays a key role in some applications in pattern recognition and machine learning;
see Devroye, Gy\"{o}rfi and Lugosi \cite{D96} and Vapnik \cite{v95}.
Denote the random VC entropy by
$$H=H(\xi_1,..., \xi_n ).$$
Boucheron,  Lugosi and Massart \cite{BLM00}
have obtained the following concentration inequalities: for any positive $t$,
\begin{eqnarray}\label{sgf01}
 \mathbf{P}\Big(\,H   -\mathbf{E}[H ]\geq t    \Big)   \leq  \exp\bigg\{ - \frac{t^2}{2 \big( \mathbf{E}[H  ]  +  \frac{1}{3}  t   \big)} \bigg\}
\end{eqnarray}
and
\begin{eqnarray}\label{sgf02}
 \mathbf{P}\Big(\,H   -\mathbf{E}[H ] \leq - t   \Big)   \leq   \exp\bigg\{ - \frac{t^2}{2 \mathbf{E}[H ] } \bigg\}.
\end{eqnarray}
Moreover, they also proved  that, for any $x \in \mathcal{X}^n,$
\begin{eqnarray*}
\sum_{i=1}^{n}\Big( H(x)-H(x^{(i)})\Big) \leq H(x);
\end{eqnarray*}
see Lemma 1 of \cite{BLM00}.  The last inequality and (\ref{coblm}) together implies that   $H$ is self-bounding.

Using (\ref{ineq47}) and (\ref{inbeg59}), we have the following result  similar to  (\ref{sgf01}) and (\ref{sgf02}).
\begin{theorem}
The random VC entropy satisfies
\begin{eqnarray}
 \emph{Var}   (H )\leq  \mathbf{E}[H],
\end{eqnarray}
and for any positive $ t  $,
\begin{eqnarray}\label{shf}
  \mathbf{P}\Big(\, H-\mathbf{E}[H]\geq t  \Big)  \leq   \exp\Bigg\{ - \frac{t^2}{2 \big(\emph{Var}(H) +  \frac{2^{3/2}}{3}  t  \sqrt{n} \big)} \Bigg\}.
\end{eqnarray}
Moreover, the last inequality holds  when replacing $H$ by $-H.$
\end{theorem}

By a simple calculation, it is easy to see that (\ref{shf}) improves the inequalities (\ref{sgf01}) and (\ref{sgf02})
for all   $t$ in a range  $0\leq t =O\big( \frac{1}{\sqrt{n} }\big(  \mathbf{E}[H]  -\textrm{Var}  (H  ) \big) \big)$ as $n\rightarrow \infty$.

\subsection{Rademacher averages}
As another application of Theorem \ref{th1ij}, we consider Rademacher averages which play an important role in the theory of
probability in Banach spaces; see Ledoux and Talagrand \cite{LTB}.
Let $B$ denote a separable Banach space, and let $X_1,...,X_n$ be independent and identically distributed bounded $B-$valued random variables.
Without loss of generality, assume that $||X_1||\leq 1$ almost surely. The
quantity of interest is the conditional Rademacher average
$$R=\mathbf{E}\bigg[ \Big|\Big|\sum_{i=1}^n\varepsilon_i X_i\Big|\Big| \bigg|  X_1,...,X_n \bigg],$$
where the $\varepsilon_i$ are independent Rademacher random variables.

 Boucheron,  Lugosi and Massart \cite{BLM03}
proved that the inequalities  (\ref{sgf01}) and (\ref{sgf02})   hold  when replacing $H$ by $R.$   (cf.\ Theorem 16).
They also proved that  $f(x)=\mathbf{E}\big[ \big|\big|\sum_{i=1}^n\varepsilon_i x_i\big|\big| \big]$ is   self-bounding.
Using (\ref{inbeg59}) and (\ref{ineq47}) again, we obtain the following concentration inequalities for $R,$ which refines the inequalities of Boucheron,  Lugosi and Massart \cite{BLM03}
for all $t$ in the range  $0\leq t =O\big( \frac{1}{\sqrt{n} }(  \mathbf{E}[R]  -\textrm{Var}  (R  ) ) \big).$
\begin{theorem}\label{th16}
The conditional Rademacher average $R$ satisfies
\begin{eqnarray}\label{sf6fd}
 \emph{Var}   (R )\leq  \mathbf{E}[R]
\end{eqnarray}
and, for any positive $ t  $,
\begin{eqnarray}\label{shfs}
  \mathbf{P}\Big(\, R-\mathbf{E}[R]\geq t  \Big)  \leq   \exp\Bigg\{ - \frac{t^2}{2 \big(\emph{Var}(R) +  \frac{2^{3/2}}{3}  t  \sqrt{n} \big)} \Bigg\}.
\end{eqnarray}
Moreover, the last inequality holds  when replacing $R$ by $-R.$
\end{theorem}

\section{Counting small subgraphs in random graphs}\label{sec5}
Consider the Erd\"{o}s-R\'{e}nyi $G(n, p)$ model of a random graph. Such
a graph has $n$ vertices and for each pair $(u, v)$ of vertices an edge is inserted between $u$ and
$v$ with probability $p$, independently. We write $m=\binom{n}{2},$  and denote the indicator variables of the
$m$ edges by $Y_1,..., Y_m.$

In this section we consider the number of triangles in a random graph. A triangle is a set of three edges
defined by vertices $u, v, w$ such that the edges are  of the form $\{u, v\},$ $\{v, w\}$ and $\{w, u\}.$
Let $Z$ denote the number of triangles in a random graph. Thus $Z$ can be expressed in the following form
$$Z= \sum_{ (i, j, k) \in \ell } Y_iY_jY_k,$$
where $\ell$ contains all triples of edges which form a triangle.
 Note  that
 $$\mathbf{E}[Z]=\binom{n}{3} p^3 \sim \frac{n^3p^3}{6} $$
and
\begin{eqnarray*}
\textrm{Var}(Z)  =   \binom{n}{3}(p^3-p^6) + 2\binom{n}{4}\binom{4}{2}(p^5-p^6)  \sim \frac{n^3}{6}(p^3-p^6) + \frac{n^4}{2}(p^5-p^6)
\end{eqnarray*}
as $n\rightarrow \infty.$
Boucheron, Lugosi and Massart \cite{BLM03} offered the following exponential inequality for
the upper tail probabilities of $Z$. Let $K>1.$ Then, for all $0\leq t \leq (K^2-1)\mathbf{E}[Z],$
\begin{eqnarray}
& &\!\!\!\!\!\!\!\! \mathbf{P}\Big(  Z- \mathbf{E}[Z]  \geq t \Big) \nonumber \\
 & &\!\!\!\!\!\!\!\!   \leq  \exp\bigg\{ -\frac{ t^2}{ (K+1)^2 \mathbf{E}[Z] (24np^2 + 24 \log n + 14t/((K+1)\sqrt{\mathbf{E}[Z]}))  }   \vee  \frac{ t^2}{ 12n \mathbf{E}[Z] +6nt  } \bigg\}. \ \ \ \ \ \ \label{BLM0}
\end{eqnarray}

We offer the following exponential inequality for tail probabilities of $Z$.
\begin{theorem} \label{thERq}
Let $Z$ denote the number of triangles in the random Erd\"{o}s-R\'{e}nyi graph.   Then,  for any positive $ t$,
\begin{eqnarray}
  \mathbf{P}\Big(\, Z -\mathbf{E}[Z]\geq t  \Big)
   \leq \exp\bigg\{ - \frac{t^2}{2 \big( \emph{Var}(Z) +  \frac{2^{3/2}}{3}   (n-2)\sqrt{n} \, t   \big)}    \bigg\}. \label{fdtyb}
\end{eqnarray}
Moreover, the same inequality holds when replacing $Z$ by $-Z.$
\end{theorem}
\emph{Proof.}  Denote  $Z =Z(Y_1,...,Y_n).$
It is easy to see that
\begin{eqnarray*}
  Z (Y_1,...,Y_k,..., Y_n) - Z (Y_1,...,Y_k',..., Y_n)
 =    \sum_{ (i, j, k) \in \ell } Y_iY_j(Y_k-Y_k').
\end{eqnarray*}
Denote by $Y_k$ the edge of the form  $\{u, v\}.$ Then $w$ belongs to the rest $n-2$ vertices, and
$ \sum_{ (i, j, k) \in \ell } Y_iY_j(Y_k-Y_k')$ contains $n-2$ summands. Since $|Y_iY_j(Y_k-Y_k')|\leq 1,$
we deduce that
\begin{eqnarray}
\sup_{k  } \Big| Z (Y_1,...,Y_k,..., Y_n) - Z (Y_1,...,Y_k',..., Y_n)\Big| \leq  n-2  .
\end{eqnarray}
Applying Theorem   \ref{co3}    to $Z,$ we obtain (\ref{fdtyb}).
  \hfill\qed

To understand inequality (\ref{fdtyb}), we summarize some of its consequences for different choices of $t$ and
for different ranges of the parameter $p$. For different ranges of $t$, we obtain the following bounds:
For all $0\leq t =o( n^{3/2}p^3 + n^{5/2}p^5 ),$
\begin{eqnarray} \label{sdsdq}
  \mathbf{P}\Big(\, Z -\mathbf{E}[Z]\geq t  \Big)
   \leq  \exp\bigg\{ - \frac{t^2}{2 \,  \textrm{Var} (Z) }  \bigg\}.
\end{eqnarray}
This is the ``Gaussian" range. Notice that the denominator coincides with the variance.
For all $ n^{3/2}p^3 + n^{5/2}p^5 \leq t  ,$
\begin{eqnarray}
  \mathbf{P}\Big(\, Z -\mathbf{E}[Z]\geq t  \Big)
   \leq \exp\bigg\{ - \frac{3 t }{2^{7/2} (n-2)\sqrt{n}  }  \bigg\}.
\end{eqnarray}

Comparing to  inequality (\ref{BLM0}), our inequality has the following three features: First,
our inequality holds for all $t\geq 0$ instead of the range $0\leq t \leq (K^2-1)\mathbf{E}[Z].$ Second, our
inequality  (\ref{sdsdq}) gives a   sharper bound. Indeed, the bound (\ref{sdsdq}) behaviors as $\exp\Big\{ - t^2/ \big(  n^4p^5+ \frac13 n^3p^3  \big)   \Big\},$ while the  bound (\ref{BLM0}) looks as $\exp\Big\{ - t^2/ \big( 4(K+1)^2n^4p^5 \wedge  2n^4p^3   \big)   \Big\}.$
Thus the bound (\ref{sdsdq}) is less than the bound  (\ref{BLM0}) for all $0\leq t =o( n^{3/2}p^3 + n^{5/2}p^5 ).$
Third, we give an upper bound for $-Z,$ that is  an upper bound on the tail probabilities $ \mathbf{P}\big(  Z- \mathbf{E}[Z]  \leq - t \big)$ for all $t>0$, while inequality  (\ref{BLM0}) usually does not provide a  similar inequality for $-Z.$

\section{Concentration for  maxima of empirical processes}\label{sec7}
\subsection{Talagrand's inequality}
Talagrand (cf.\ Theorem 1.4 of \cite{T96}) gave the following concentration inequalities for the maxima of empirical processes.
\begin{theorem}\label{sgff}
  Let $(\xi_i)_{i=1,...,n}$ be a sequence of independent random variables with values
in a measurable space  $\mathcal{X}.$ Consider a countable class $\mathcal{F}$ of measurable
functions on $\mathcal{X}.$ Consider the random variable $$Z=\sup_{f \in \mathcal{F}} \sum_{i=1}^n f(\xi_i).$$
Denote
\[
U= \sup_{f \in \mathcal{F}} ||f(\xi_i)||_{\infty}\ \ \  \ \ \textrm{and} \ \ \ \ \ V= \mathbf{E}\Big[\sup_{f \in \mathcal{F}} \sum_{i=1}^n f^2(\xi_i)\Big].
\]
Then,  for any  positive $t,$
\begin{eqnarray} \label{TLAin}
\mathbf{P}\Big(  Z- \mathbf{E}[Z]  \geq t \Big) \, \leq \,  K \exp\bigg\{ -\frac{t^2}{2(c_1V+ c_2 U t ) }   \bigg\},
\end{eqnarray}
where $K,   c_1$ and $c_2$ are positive absolute constants.  Moreover, the same inequality  holds when replacing $Z$ by $-Z.$
\end{theorem}

 Talagrand's proof of Theorem \ref{sgff} is rather intricate and does not lead to very attractive values for the constants $K,   c_1$ and $c_2.$
 One year after Talagrand's work \cite{T96},
Ledoux \cite{L96}   developed a new and much simpler method for establishing similar concentration inequalities of Talagrand.
Ledoux's method is known as ``entropy method''.
His method allows one to obtain  explicit constants.  In particular, Ledoux \cite{L96} showed that for an adequate positive constant $C$, replacing $V$ by $V + C U  \mathbf{E}[Z],$
inequality (\ref{TLAin}) holds with $K=2, c_1=42$ and $c_2=8.$  Notice that, by a remark of  Massart \cite{M00}, Ledoux's inequalities did not recover exactly Talagrand's inequality (\ref{TLAin}) due to the difference of $V$.   Moreover, he did not provide the same bound for $-Z$ as Talagrand's inequalities  in general.
The entropy method has become very popular in recent years. Many interesting concentration inequalities have been established via this method; see Bobkov and Ledoux \cite{BL97}, Massart \cite{M00},    Rio \cite{R01,R02}, Bousquet \cite{B02},   Boucheron \emph{et al.} \cite{BLM03},  Klein and Rio \cite{KR05}.
Similar to Ledoux's inequality, such type concentration inequalities also do not hold for $Z$ replacing by $-Z$ in general.

\subsection{Talagrand type concentration  inequalities}
In this section, we would like to give some concentration  inequalities for the maxima of empirical processes. Our inequalities are similar to Talagrand's inequality (\ref{TLAin}), and they hold also when replacing $Z$ by $-Z.$ Moreover,   our inequalities are  Gaussian   bounds. To the best of our knowledge, such type bounds for the maxima of empirical processes have not been obtained before.

\begin{theorem} \label{supeinq}
  Let $(\xi_i)_{i=1,...,n}$ be a sequence of independent random variables with values
in a measurable space  $\mathcal{X}.$ Consider a countable class $\mathcal{F}$ of measurable
functions on $\mathcal{X}.$ Assume that
\[
U_i= \sup_{f \in \mathcal{F}} ||f(\xi_i)||_{\infty} \ < \infty.
\]
Let $Z$ denote as one of the following formulas
$$ \sup_{f \in \mathcal{F}} \sum_{i=1}^n f(\xi_i),  \ \ \ \ \   \ \ \ \ \  \sup_{f \in \mathcal{F}} \Big|\sum_{i=1}^n f(\xi_i)\Big|,$$
$$   \sup_{f \in \mathcal{F}}  \sum_{i=1}^n f(\xi_i) - \mathbf{E}  \big[   f(\xi_i) \big] \ \ \ \ \ \ \textrm{and} \  \ \ \ \ \ \ \sup_{f \in \mathcal{F}} \Big|\sum_{i=1}^n  f(\xi_i)  - \mathbf{E}\big[   f(\xi_i)   \big]\Big|    .$$
Denote   $$ U= \max\{ U_i:\, i=1,...,n\} \ \  \ \  \ \textrm{and} \   \ \ \    \ \sigma^2= \frac1n \emph{Var}(Z).$$
Then, for  any positive $t,$
\begin{eqnarray}
 \mathbf{P}\Big(  Z- \mathbf{E}[Z]  \geq t \sqrt{\emph{Var}(Z)} \Big)
&\leq&  \exp\Bigg\{ - \frac{t^2}{ 1+ \sqrt{1+\frac{2^{7/2}  }{3 } U t \big/\sigma  } \, + \frac{ 2^{5/2} }{ 3  } U t/\sigma   } \Bigg\} \ \ \ \ \label{fdfsda}\\
 &\leq&  \exp\bigg\{ -\frac{t^2}{2 \big( 1 +     \frac{2^{5/2}}3 U t/ \sigma  \big) }   \bigg\}.\label{fdfsdb}
\end{eqnarray}
Moreover, the same inequalities hold when replacing $Z$ by $-Z.$
\end{theorem}
\emph{Proof}.  Denote by $Z =Z_n(\xi_1,...,\xi_n).$ When $Z  = \sup_{f \in \mathcal{F}} \sum_{i=1}^n f(\xi_i),$ by the fact $U_i= \sup_{f \in \mathcal{F}} ||f(\xi_i)||_{\infty},$ it is easy to see that
\begin{eqnarray*}
  Z_n(\xi_1,...,\xi_k,..., \xi_n) - Z_n(\xi_1,...,\xi_k',..., \xi_n)
&\geq& Z  - \sup_{f \in \mathcal{F}} \Big(\sum_{i=  1}^n  f(\xi_i)+  f(\xi_k')-f(\xi_k) \Big) \\
&\geq& Z  - Z - \sup_{f \in \mathcal{F}}  \Big( f(\xi_k')-f(\xi_k) \Big) \\
&\geq&   -  2\, U_k,
\end{eqnarray*}
where $(\xi_k')_{k=1,...,n}$ is an independent copy of $(\xi_k)_{k=1,...,n}.$
The last inequality holds also for $$ Z_n(\xi_1,...,\xi_k',..., \xi_n)- Z_n(\xi_1,...,\xi_k,..., \xi_n).$$
Thus
\begin{eqnarray}\label{uneg53}
\sup_{\stackrel{ \xi_1,...,\xi_n}{ \xi_1',...,\xi_n'\in \mathbf{E}}     } \Big|Z_n(\xi_1,...,\xi_k,..., \xi_n) - Z_n(\xi_1,...,\xi_k',..., \xi_n)\Big| \leq 2 \, U_k.
\end{eqnarray}
Applying Theorem \ref{co1}  to $Z,$ we obtain (\ref{fdfsda}) and (\ref{fdfsdb}) for $Z = \sup_{f \in \mathcal{F}} \sum_{i=1}^n f(\xi_i).$  When  $Z  = \sup_{f \in \mathcal{F}} |\sum_{i=1}^n f(\xi_i)|,$ by the fact $\sup_{f \in \mathcal{F}}   ||f(\xi_k') - f(\xi_k)||_{\infty} \leq 2\, U_k $ again,
the inequalities (\ref{fdfsda}) and (\ref{fdfsdb}) also follow by a similar argument.
  When $Z$ denotes either
$$   \sup_{f \in \mathcal{F}}  \sum_{i=1}^n f(\xi_i) - \mathbf{E}  \big[   f(\xi_i) \big] \ \ \ \ \ \ \textrm{or} \  \ \ \ \ \ \ \sup_{f \in \mathcal{F}} \Big|\sum_{i=1}^n  f(\xi_i)  - \mathbf{E}\big[   f(\xi_i)   \big]\Big| ,$$
by the fact
\[
 \sup_{f \in \mathcal{F}} \Big|\Big|f (\xi_i) - \mathbf{E}  \big[   f(\xi_i ) \big] - f (\xi_i') + \mathbf{E}  \big[   f(\xi_i') \big]\Big|\Big|_{\infty}= \sup_{f \in \mathcal{F}} \Big|\Big|f (\xi_i) -    f(\xi_i')  \Big|\Big|_{\infty}  \leq  2\, U_i,
\]
inequality (\ref{uneg53}) holds true.  Then the inequalities (\ref{fdfsda}) and (\ref{fdfsdb}) follow  by a  similar argument again.      \hfill\qed

\begin{remark}
Inspiring the proof of Theorem \ref{supeinq}, it is easy to see that  if denote $Z$ as one of the following formulas
$$ \sup_{(f_1,...,f_n) \in \mathcal{F}^n} \sum_{i=1}^n f_i(\xi_i),  \ \ \ \ \   \ \ \ \ \  \sup_{(f_1,...,f_n) \in \mathcal{F}^n} \Big|\sum_{i=1}^n f_i(\xi_i)\Big|,$$
$$   \sup_{(f_1,...,f_n) \in \mathcal{F}^n}  \sum_{i=1}^n f_i(\xi_i) - \mathbf{E}  \big[   f_i(\xi_i) \big] \ \ \ \ \ \ \textrm{and} \  \ \ \ \ \ \ \sup_{(f_1,...,f_n) \in \mathcal{F}^n} \Big|\sum_{i=1}^n  f_i(\xi_i)  - \mathbf{E}\big[   f_i(\xi_i)   \big]\Big| ,$$
then the inequalities
(\ref{fdfsda}) and  (\ref{fdfsdb}) also hold true.
\end{remark}

Notice that (\ref{fdfsdb}) can  be rewritten in the following form: for any positive $x$,
 \begin{eqnarray}\label{ineq65}
\mathbf{P}\Big(  Z- \mathbf{E}[Z]  \geq x  \Big)   \leq    \exp\bigg\{ -\frac{x^2}{2 \big(   \textrm{Var}   (Z) +     \frac{2^{5/2}}3 U x \sqrt{n} \big) }   \bigg\}.
\end{eqnarray}
The last inequality shows that we can take $K=1, c_1=1,$ $V= \textrm{Var}   (Z)$ and $c_2 =\frac{2^{5/2}}3  \sqrt{n}$ in Talagrand's inequality (\ref{TLAin}).
This provides a partial positive  answer to Massart's question \cite{M00} about the best constants in   Talagrand's  inequality.
Moreover, our bound (\ref{ineq65}) is less than Talagrand's bound (\ref{TLAin}) for $x$ in the range  $$0 \leq x =O\bigg(\frac{ 1}{  U\sqrt{n}}\Big( \mathbf{E}\Big[\sup_{f \in \mathcal{F}} \sum_{i=1}^n f^2(\xi_i)\Big] - \textrm{Var} (Z)\Big) \bigg).$$
When $\textrm{Var} (Z) / \mathbf{E}\big[\sup_{f \in \mathcal{F}} \sum_{i=1}^n f^2(\xi_i)\big] \rightarrow 0$ as $n\rightarrow \infty, $ this range is large.  To illustrate it, consider the case that $\textrm{Var} (Z)$ is in order of $n^\alpha, 0< \alpha <1,$ and $\mathbf{E}\big[\sup_{f \in \mathcal{F}} \sum_{i=1}^n f^2(\xi_i)\big]$ is in order of $n$ as $n\rightarrow \infty,$ then our inequalities is significantly better than Talagrand's inequality for \emph{standard} $t$ in the range $0\leq t =o(n^{(1-\alpha)/2})$.
Recall that a $t$ is  called standard if we refer to the tail probabilities $\mathbf{P}\big( ( Z- \mathbf{E}[Z]) /\sqrt{\textrm{Var} (Z)} \geq t \big).$

Here we would like to give some  comparisons between our results and the inequalities of  Klein and Rio \cite{KR05}.  Assume that $ \mathbf{E}[f(\xi_i)]=0$ for all $i\in[1, n]$ and all $f \in \mathcal{F}.$ Klein and Rio \cite{KR05} have obtained the following inequality.
  Let $Z=\sup_{f \in \mathcal{F}} \sum_{i=1}^n f(\xi_i)$ and let $U=1$.
Then, for  any positive $x,$
\begin{eqnarray} \label{sdgfmff}
\mathbf{P}\Big(  Z- \mathbf{E}[Z]  \geq x \Big) \, \leq \,    \exp\bigg\{ -\frac{x^2}{2 v +     3 \, x  }   \bigg\}
\end{eqnarray}
and
\begin{eqnarray} \label{sdgfmff02}
\mathbf{P}\Big(   Z - \mathbf{E}[Z]  \leq - x \Big) \, \leq \,    \exp\bigg\{ -\frac{x^2}{2 v +     2 \, x   }   \bigg\},
\end{eqnarray}
where  $$v= \sup_{f \in \mathcal{F}}    \textrm{Var}   \Big(\sum_{i=1}^n f (\xi_i)\Big) + 2\,\mathbf{E}[Z]  $$
and $\sup_{f \in \mathcal{F}}    \textrm{Var}   \big(\sum_{i=1}^n f (\xi_i)\big) \leq\mathbf{E}\big[\sup_{f \in \mathcal{F}} \sum_{i=1}^n f^2(\xi_i)\big].$
Moreover, Klein and Rio (cf.\ Corollary 1.1 of \cite{KR05}) have pointed out that
\begin{eqnarray} \label{sdgfmff03}
 \textrm{Var} (Z) \leq v.
\end{eqnarray}
Comparing to the results of Klein and Rio \cite{KR05}, our inequalities have the following three interesting features.
 First, our inequalities (\ref{fdfsda}) and (\ref{fdfsdb}) do not need the assumption that $(\xi_i)_{i=1,...,n}$  are  centered with respect to $f,$ i.e.\ $ \mathbf{E}[f(\xi_i)]=0$ for all $i\in[1, n]$ and all $f \in \mathcal{F}.$ Second,  it is easy to see that our inequalities (\ref{fdfsda}) and (\ref{fdfsdb})  are  Gaussian type bounds, while the inequality of Klein and Rio \cite{KR05} does not share this feature due to the fact $\textrm{Var} (Z) \leq v.$
Third, by (\ref{ineq65}), we find  that our bound is better than Klein-Rio's bound  for $x$ in the range  $0 \leq x =\frac{ 3}{2^{5/2} \sqrt{n}}\big( v - \textrm{Var} (Z)\big).$ For instance, if $v$ is in order of $n$ and $\textrm{Var} (Z)$ is in order of $n^\alpha, 0< \alpha <1,$ as $n\rightarrow \infty,$ then for \emph{standard} $t$ in the range $0\leq t =o(n^{(1-\alpha)/2})$ our bound is better than Klein-Rio's.

Inequality (\ref{fdfsdb})   can also be rewritten in the following form: for any positive $t,$
\begin{eqnarray} \label{fdsclp}
\mathbf{P}\Big(   Z- \mathbf{E}[Z]   \geq \sqrt{n} \big(  \sqrt{c^2U^2 t^2 +2 t   \sigma^2 }+ c \, U t \big)\Big) \, \leq \,    \exp\Big\{ -t \Big\},
\end{eqnarray}
where $c=2^{5/2}/3.$
Since $$\sqrt{c^2U^2 t^2 +2 t  \sigma^2 } \leq  \sigma \sqrt{ 2 \, t  }  + c \, U  t ,$$
inequality (\ref{fdsclp}) implies the following bound: for any positive $t,$
\begin{eqnarray} \label{fdgtsclp}
\mathbf{P}\Big(   Z- \mathbf{E}[Z]    \geq  \sqrt{n} \big(\sigma \sqrt{ 2 \, t  }  + 2\, c\, U t  \big)  \Big) \, \leq \,    \exp\Big\{ -t \Big\}.
\end{eqnarray}
Such type  bound can be found in  Theorem 3  of Massart \cite{M00}.
When $Z$ denotes either
$$  \sup_{f \in \mathcal{F}} \Big|\sum_{i=1}^n f(\xi_i)\Big| \ \ \   \ \textrm{or} \   \ \ \ \sup_{f \in \mathcal{F}} \Bigg|\sum_{i=1}^n  f(\xi_i)  - \mathbf{E}\big[   f(\xi_i)   \big]\Bigg|  ,$$
 Massart has established the following inequality: for positive $t $ and a positive $\varepsilon,$
\begin{eqnarray}
\mathbf{P}\Big(   Z- \mathbf{E}[Z]    \geq  \varepsilon \, \mathbf{E}[Z] + \varrho \sqrt{ 2 \,\kappa \, t  }  + \kappa(\varepsilon) \, U t    \Big) \, \leq \,    \exp\Big\{ -t \Big\}, \nonumber
\end{eqnarray}
where $\varrho^2 = \sup_{f \in \mathcal{F}}  \sum_{i=1}^{n}   \textrm{Var} (f (\xi_i)), \kappa=4$ and $\kappa(\varepsilon)=  2.5+ 32\, \varepsilon^{-1}.$


\section*{Acknowledgements}
Fan would like to thank the anonymous referee and editor for their helpful comments and suggestions.

\end{document}